\documentstyle{amsppt}
\NoBlackBoxes

\define\cc{\Bbb C}
\define\cn{\Bbb C^n}
\define\cm{\Bbb C^m}
\define\cd{\Bbb C^2}

\define\er{\Bbb R}

\define\wl{\Cal L_\infty}

\define\la{\lambda}
\define\vf{\varphi}

\define\oo#1{\overset {{}_o}\to {#1}}

\document
\topmatter 
\footnote"{}"{1991 {\it Mathematics Subject Classification}: 14E05.} 
\footnote"{}"{{\it Key words and phrases}: polynomial mapping, \L ojasiewicz
exponent.} 
\footnotetext"{}"{This research was partially supported by KBN Grant No. 2
P03A 050 10}
\endtopmatter

\vglue 1.5cm
\centerline{\bf A set on which the \L ojasiewicz exponent}

\centerline{\bf at infinity is attained}

\vskip.5cm
\font\tw=cmcsc10
\centerline{\tw {\rm by} Jacek Ch\c adzy\'nski and Tadeusz Krasi\'nski {\rm (\L
\'od\'z)}} 


\font\bbb=cmr8
\font\ccc=cmbx8

\vskip.7cm
\vbox{{\ccc Abstract. }\bbb We show that for a polynomial mapping
$F=(f_1,...,f_m): \cn\to\cm$ the \L ojasiewicz exponent $\wl(F)$ of $F$ is
attained on the set $\{z\in \cn:\,f_1(z)\cdot \ldots \cdot f_m(z)=0\}$.}

\vskip.5cm 
{\bf 1. Introduction.} The purpose of this paper is to prove that the \L
ojasiewicz exponent at infinity of a polynomial mapping $F:\cn\to\cm$ is
attained on a proper algebraic subset of $\cn$ defined by the components of
$F$ (Thm 1). As a corollary we obtain a result of \/ Z. Jelonek on testing sets
for properness of polynomial mappings (Cor.\hskip 1mm 3) and a formula for the \L
ojasiewicz exponent at infinity of $F$ in the case $n=2$, $m\ge 2$, in terms
of parametrizations of branches (at infinity) of zeroes of the components of
$F$ (Thm 2). This result is a generalization of the authors' result for
$n=m=2$ ([CK], Main Thm).

Before the main considerations we show some basic properties of the \L
ojasiewicz exponent at infinity for regular mappings i.e.
for polynomial mappings restricted to algebraic subsets of $\cn$. We prove
that the exponent is a rational number, that it is attained on a meromorphic
curves (Prop. 1) and we give a condition equivalent to the properness of 
regular mappings (Cor. 2). These properties are analogous to the ones, known
in folklore, for polynomial mappings from $\cn$ into $\cm$. We do not pretend
to the originality of proof methods of these properties but we only want
to fill gaps in the literature.

The results obtained by Z. Jelonek in [J] have played the inspiring role in
the undertaking this research. On the other hand, the idea of the proof of
the main theorem was taken from A. P\l oski ([P$_2$], App.).

\vskip.3cm
{\bf 2. The \L ojasiewicz exponent.} Let $F:\cn\to\cm$, $n\ge 2$, be a
polynomial mapping and let $S\subset \cn$ be an unbounded algebraic set. Put
$$
N(F|S):= \{\nu\in \er :\exists A>0,\,\exists B>0,\,\forall z\in
S,\,(|z|>B\Rightarrow A|z|^\nu \le |F(z)|)\},
$$
where $|\cdot |$ is the polycylindric norm. If $S=\cn$ we define
$N(F):=N(F|\cn)$. 

By the \L ojasiewicz exponent at infinity of $F|S$ we mean $\wl(F|S):=\sup
N(F|S)$. Analogously $\wl(F):=\sup N(F)$.

Before we pass to properties of the \L ojasiewicz exponent we quote the known
curve selection lemma at infinity (cf [NZ], Lemma 2). We begin with a
definition. A curve $\vf:(R,+\infty)\to \er ^k$ is called meromorphic at
$+\infty$ if $\vf$ is the sum of a Laurent series of the form
$$
\vf(t)=\alpha_pt^p+\alpha_{p-1}t^{p-1}+\ldots,\qquad \alpha_i\in \er^k.
$$
By $\|\cdot\|$ we denote the euclidian norm in $\er^k$.

\proclaim{\indent{\tw Lemma 1 (Curve Selection Lemma)}} If $X\subset \er^k$ is an unbounded
semi-algebraic set, then there exists a curve $\vf:(R,+\infty)\to \er^k$,
meromorphic at $+\infty$, such that $\vf(t)\in X$ for $t\in(R,+\infty)$ and
$\|\vf(t)\|\to +\infty$ when $t\to +\infty$.
\endproclaim

Let us notice that the \L ojasiewicz exponent at infinity of a regular
mapping $F|S$ does not depend on the norm in $\cn$. So, in the sequel of this
section, we shall use the euclidian norm $\|\cdot \|$ in the definition of
$N(F|S)$. 

Let us introduce one more definition. A curve
$\vf=(\vf_1,...,\vf_m):\{t\in\cc: |t|>R\}\to \cm$ is called meromorphic at
$\infty$, if $\vf_i$ are meromorphic at $\infty$.

Let $F:\cn\to \cm$, $n\ge 2$, be a polynomial mapping and let $S\subset \cn$
be an unbounded algebraic set.

\proclaim{\indent{\tw Proposition 1}} If \; $\#(F|S)^{-1}(0)<\infty$, then
$\wl(F|S)\in N(F|S)\cap \Bbb Q$. Moreover, there exists a curve $\vf:\{t\in
\cc:|t|>R\}\to \cm$, meromorphic at $\infty$, such that $\vf(t)\in S$,
$\|\vf(t)\|\to +\infty$ for $t\to \infty$ and
$$
\|F\circ \vf(t)\|\sim \|\vf(t)\|^{\wl(F|S)}\qquad\text{for }t\to \infty.
\tag 1
$$
\endproclaim

\demo{Proof} Let us notice first that the set
$$
\{(z,w)\in S\times S:\|F(z)\|^2\le \|F(w)\|^2\; \vee\; \|z\|^2\ne \|w\|^2\}
$$
is semi-algebraic in $\cn\times \cn\cong \er^{4n}$. Then by the
Tarski-Seidenberg theorem (cf [BR], Rem. 3.8) the set
$$
\align
X:&=\{z\in S:\forall w\in S,\,(\|F(z)\|^2\le\|F(w)\|^2\; \vee \;
\|z\|^2\ne\|w\|^2)\}\\ 
&=\{z\in S:\|F(z)\|=\min_{\|w\|=\|z\|}\|F(w)\|\}
\endalign
$$
is also semi-algebraic and obviously unbounded in $\cn\cong \er^{2n}$. So,
by Lemma 1 there exists a curve $\tilde \vf:(R, +\infty)\to X$, meromorphic
at $\infty$, such that $\|\tilde \vf(t)\|\to +\infty$ for $t\to +\infty$. Then
there exists a positive integer $p$ such that $\tilde \vf$ is the sum of a
Laurent series
$$
\tilde\vf(t)=\alpha_pt^p+\alpha_{p-1}t^{p-1}+\ldots ,\qquad \alpha_i\in
\cn,\quad\alpha_p\ne0. 
\tag 2
$$
Since $\#(F|S)^{-1}(0)<\infty$ then there exists an integer $q$ such that
$F\circ \tilde \vf$ is the sum of a Laurent series
$$
F\circ \tilde \vf(t)=\beta_qt^q+\beta_{q-1}t^{q-1}+\ldots,\qquad
\beta_i\in\cm,\quad \beta_q\ne 0.
\tag 3
$$
From (2) and (3) we have
$$
\|F\circ \tilde \vf(t)\|\sim \|\tilde \vf(t)\|^\la\qquad\text{for
}t\to+\infty, 
\tag 4
$$
where $\la:=q\slash p$. Let $\tilde \Gamma:=\{z\in \cn:z=\tilde \vf(t),\,t\in
(R,+\infty)\}$. Then from (4)
$$
\|F(z)\|\sim \|z\|^\la\qquad \text{for }\|z\|\to \infty,\quad z\in
\tilde\Gamma. 
\tag 5
$$

Now, we shall show that $\wl(F|S)=\la$. From (5) we have $\wl(F|S)\le \la$.
Since $\tilde\Gamma\subset X$ is unbounded, then there exist positive constants $A$, $B$
such that $\|F(z)\|\ge A\|z\|^\la$ for every $z\in S$ and $\|z\|>B$. Then
$\la\in N(F|S)$ and in consequence $\wl(F|S)\ge \la$. Summing up,
$\wl(F|S)=\la\in N(F|S)\cap \Bbb Q$. 

Now, we shall prove the second part of
the assertion. Let $\vf$ be an extension of $\tilde \vf$ to the complex
domain, that is
$$
\vf(t)=\alpha _pt^p+\alpha_{p-1}t^{p-1}+\ldots ,
\tag 6
$$
where $t\in\cc$ and $|t|>R$. Obviously, series (6) is convergent and, as
above, $\alpha _i\in\cn$, $\alpha_p\ne 0$. Hence $\vf$ is a curve,
meromorphic at $\infty$, and clearly $\|\vf(t)\|\to+\infty$ for $t\to
\infty$. Moreover, $F\circ \vf$ is an extension of $F\circ \tilde \vf$ to the
complex domain and
$$
F\circ \vf(t)=\beta_qt^q+\beta_{q-1}t^{q-1}+\ldots ,
\tag 7
$$
where $t\in \cc$ and $|t|>R$. Obviously, the series (7) is convergent and, as
above, $\beta_i\in\cm$, $\beta_q\ne 0$. From (6), (7) and the definition of
$\la$ we get (1). Since $S$ is an algebraic subset of $\cn$ and $\tilde
\vf(t)\in S$ for $t\in (R,+\infty)$, then also $\vf(t)\in S$ for $t\in \cc$,
$|t|>R$.

This ends the proof of the proposition.
\enddemo

Let $F:\cn\to\cm$, $n\ge 2$, be a polynomial mapping and $S\subset \cn$ -- an
algebraic unbounded set.

Directly from Proposition 1 we get

\proclaim{\indent{\tw Corollary 1}} $\wl(F|S)>-\infty$ if and only if
$\#(F|S)^{-1}(0)<\infty$. 
\endproclaim

We get also easily from Proposition 1

\proclaim{\indent{\tw Corollary 2}} The mapping $F|S$ is proper if and only if
$\wl(F|S)>0$.
\endproclaim

In fact, if $\wl(F|S)>0$, then obviously $F|S$ is a proper mapping. If, in
turn, $\wl(F|S)\le 0$ then from the second part of Proposition 1 and
Corollary 1 it follows that there exists a sequence $z_n\in S$ such that
$\|z_n\|\to +\infty$ and the sequence $F(z_n)$ is bounded. Hence $F|S$ is not
proper mapping in this case.

\vskip.3cm
{\bf 3. The main result.} Now, we formulate the main result of the paper.

\proclaim{\indent{\tw Theorem 1}} Let $F=(f_1,...,f_m):\cn\to\cm$, $n\ge 2$
be a polynomial mapping and $S:=\{z\in \cn:\,f_1(z)\cdot\ldots\cdot
f_m(z)=0\}$. If $S\ne \emptyset$, then
$$
\wl(F)=\wl(F|S).
\tag 8
$$
\endproclaim

The proof will be given in section 4.

Directly from Theorem 1 and Corollary 2 we get

\proclaim{\indent{\tw Corollary 3} {\rm ([J], Cor. 6.7)}} If
$F=(f_1,...,f_m):\cn\to\cm$, $n\ge 2$ is a polynomial mapping and $S:=\{z\in
\cn:f_1(z)\cdot\ldots\cdot f_m(z)=0\}$ is not empty, then the mapping $F$ is
proper if and only if the mapping $F|S$ is proper.
\endproclaim

Another corollary from Theorem 1 is an effective formula for the \L
ojasiewicz exponent, generalizing the earlier authors' result ([CK], Main
Thm). 

Let us introduce some notions. If $\Psi:\{z\in \cc:|z|>R\}\to \cc^k$ is the
sum of a Laurent series of the form 
$$
\Psi(t)=\alpha_pt^p+\alpha_{p-1}t^{p-1}+\ldots ,\qquad \alpha_i\in
\cc^k,\quad \alpha_p\ne 0
$$
then we put $\deg \Psi:=p$. Additionally, $\deg \Psi:=-\infty$ if $\Psi =0$.
For an algebraic curve in $\cd$, the notions of its branches in a
neighbourhood of $\infty$ and parametrizations of these branches we take
after [CK].

Let now $F=(f_1,...,f_m):\cd\to\cm$ be a polynomial mapping and $S:=\{z\in
\cd:f_1(z)\cdot\ldots\cdot f_m(z)=0\}$. Assume that $S\ne \emptyset$ and $S\ne
\cd$. 

\proclaim{\indent{\tw Theorem 2}} If $\Gamma_1$,...,$\Gamma_s$ are branches
of the curve $S$ in a neighbourhood of infinity $Y$ and $\Phi_i:U_i\to Y$,
$i=1,...,s$, are their parametrizations, then
$$
\wl(F)=\min_{i=1}^s\frac{\deg F\circ \Phi_i}{\deg \Phi_i}.
\tag 9
$$
\endproclaim

\demo{Proof} Denote $\la_i:=\deg F\circ \Phi_i\slash \deg \Phi_i$. If
$\la_i=-\infty$ for some $i$, then (9) holds. So, assume that
$\la_i\ne-\infty$, $i=1,...,s$. Then
$$
|F(z)|\sim |z|^{\la_i}\qquad \text{for }|z|\to+\infty,\quad z\in \Gamma_i.
$$
Hence, taking into account the equality $S\cap Y=\Gamma_1\cup\ldots \cup
\Gamma_s$ we get (9).
\enddemo

\vskip.3cm
{\bf 4. Proof of the main theorem.} Let us begin with a lemma on polynomial
mappings from $\cc$ into $\cm$. It is a generalization of a result by A. P\l
oski ([P$_1$], Lemma 3.1) and plays a key role in the proof of the
main theorem.

\proclaim{\indent{\tw Lemma 2}} Let $\Phi =(\vf_1,...,\vf_m):\cc\to \cm$ be a
polynomial mapping and $\vf:=\vf_1\cdot\ldots\cdot \vf_m$. If $\vf$ is a
polynomial of positive degree and $T$ is its set of zeroes, then for every $t\in
\cc$ 
$$
|\Phi(t)|\ge 2^{-\deg \Phi}\min_{\tau\in T}|\Phi(\tau)|.
$$
\endproclaim

\demo{Proof} Fix $t_0\in \cc$. Let $\min_{\tau\in T}|t_0-\tau|$ be attained
for some $\tau_0\in T$. If $\vf_i$ is a polynomial of positive degree and has
the form $\vf_i(t)=c_i\prod^{\deg \vf_i}_{j=1}(t-\tau_{ij})$, then we have
$$
2|t_0-\tau_{ij}|=|t_0-\tau_{ij}|+|t_0-\tau_{ij}|\ge
|t_0-\tau_0|+|t_0-\tau_{ij}|\ge |\tau_0-\tau_{ij}|.
$$
Hence
$$
2^{\deg \vf_i}|\vf_i(t_0)|\ge |\vf_i(\tau_0)|.
$$
Obviously, this inequality is also true for $\vf_i$ being a constant. Since
$\deg \Phi\ge \deg \vf_i$, then from the above we get
$$
2^{\deg \Phi}|\Phi(t_0)|\ge |\Phi(\tau_0)|\ge \min_{\tau\in T}|\Phi(\tau)|,
$$
which ends the proof.
\enddemo

In the sequel, $z=(z_1,...,z_n)\in \cn$, $n\ge 2$, and for every $i\in
\{1,...,n\}$ we put $z_i':=(z_1,...,z_{i-1},z_{i+1},...,z_n)$.

Without proof let us notice an easy lemma.

\proclaim{\indent{\tw Lemma 3}} Let $f:\cn\to\cc$ be a non-constant
polynomial function and $S$ -- its set of zeroes. If for every $i\in
\{1,...,n\}$, $\deg f=\deg _{z_i}f$, then there exist constants $C\ge 1$,
$D>0$ such that for every $i\in\{1,...,n\}$,
$$
|z_i|\le C|z'_i|\qquad\text{for }z\in S\text{ and }|z_i'|>D.
$$
\endproclaim

\demo{Proof of Theorem 1} Without loss of generality we may assume that

(i) $S\ne \cn$,

\noindent and that

(ii) $\#(F|S)^{-1}(0)<\infty$.

\noindent In fact, if (i) does not hold then (8) is obvious, whereas, if (ii)
does not hold then (8) follows from Corollary 1.

Obviously $N(F)\subset N(F|S)$. So, to prove (8) it suffices to show
$$
N(F|S)\subset N(F).
\tag 10
$$

Put $f:=f_1\cdot\ldots\cdot f_m$. From (i) we have $\deg f>0$. Since the sets
$N(F|S)$ and $N(F)$ are invariant with respect to linear changes of
coordinates in $\cn$ we may assume that
$$
\deg f=\deg_{z_i} f,\qquad i=1,...,n.
\tag 11
$$
This obviously implies
$$
\deg f_j=\deg _{z_i}f_j,\qquad j=1,...,m,\quad i=1,...,n.
\tag 12
$$

It follows from (ii) and Corollary 1 that $N(F|S)$ is not empty. Take $\nu\in
N(F|S)$. Then there exist $A>0$, $B>0$ such that
$$
|F(\zeta)|\ge A|\zeta |^\nu,\qquad \text{for }\zeta \in S,\quad |\zeta|>B.
\tag 13
$$
By (11) and Lemma 3 there exist $C\ge 1$, $D>0$ such that
for every $i\in \{1,...,n\}$,
$$
|z_i|\le C|z_i'|\qquad\text{for }z\in S,\quad |z_i'|>D.
\tag 14
$$
Put $A_1:=2^{-\deg F}A\min(1,C^\nu)$\/ and\/ $B_1:=\max (B,D)$. Take arbitrary
$\oo{z}\in\cn$ such that $|\oo{z}|>B_1$. Clearly, $|\oo{z}|=|\oo{z}_i'|$ for
some $i$. Define
$\vf_j(t):=f_j(\oo{z}_1,...,\oo{z}_{i-1},t,\oo{z}_{i+1},$ $...,\oo{z}_n)$,
$\Phi:=(\vf_1,...,\vf_m)$. Then from (12) we have 
$$
\deg F=\deg \Phi.
\tag 15
$$
Moreover, from (11) it follows that $\vf:=\vf_1\cdot\ldots\cdot \vf_m$ is a
polynomial of positive degree. Then, from Lemma 2 ($T$ is defined as in Lemma
2) and (15) we have
$$
|F(\oo{z})|=|\Phi(\oo{z}_i)|\ge 2^{-\deg \Phi}\min_{\tau\in
T}|\Phi(\tau)|=2^{-\deg F}|F(\oo{\zeta})|
\tag 16
$$
for some $\oo{\zeta}=(\oo{z}_1,...,\oo{z}_{i-1},\tau_0,\oo{z}_{i+1},...,\oo{z}_n)$,
$\tau _0\in T$. So, $\oo{\zeta}\in S$. Since $|\oo{z}|>B_1$ and $|\oo{\zeta}|\ge
|\oo{z}_i'|=|\oo{z}|$, then from (16) and (13) we get
$$
|F(\oo{z})|\ge 2^{-\deg F}A|\oo{\zeta}|^\nu,
\tag 17
$$
whereas from (14)
$$
|\oo{z}|\le |\oo{\zeta}|\le C|\oo{z}|.
\tag 18
$$
Considering two cases, when $\nu\ge 0$ and $\nu<0$, from (17) and (18) we
easily get
$$
|F(\oo{z})|\ge A_1|\oo{z}|^\nu.
$$
Since $\oo{z}$ is arbitrary we have $\nu\in N(F)$.

This ends the proof of the theorem.
\enddemo

\font\bbb=cmr8

\vskip.3cm
\noindent{\baselineskip=10pt
\bbb Faculty of Mathematics\newline
University of \L \'od\'z\newline
S. Banacha 22\newline
90-238 \L \'od\'z, Poland

\noindent
E-mail: jachadzy\@imul.uni.lodz.pl

\noindent\hskip 11mm krasinsk\@krysia.uni.lodz.pl}

\enddocument